\documentclass[12pt,leqno]{article}
\usepackage{enumerate}
\usepackage{a4wide}
\usepackage{amsmath}
\usepackage{amssymb}
\usepackage[all]{xy}
\usepackage[T1]{fontenc}
\usepackage[latin1]{inputenc}

\newtheorem{thm}    {Theorem}
\newtheorem{prop}   [thm]{Proposition}
\newtheorem{lem}    [thm]{Lemma}
\newtheorem{cor}    [thm]{Corollary}
\newtheorem{defi}   [thm]{Definition}

\newcommand{\CC}{\mathbb{C}}
\newcommand{\QQ}{\mathbb{Q}}
\newcommand{\FF}{\mathbb{F}}
\newcommand{\ZZ}{\mathbb{Z}}

\newcommand{\Qbar}{{\overline{\QQ}}}

\newcommand{\Fpbar}{{\overline{\FF_p}}}
\newcommand{\Ftbar}{{\overline{\FF_2}}}
\newcommand{\Katz}{\mathrm{Katz}}

\newcommand{\univ}{\mathrm{univ}}
\newcommand{\Norm}{\mathrm{Norm}}

\newcommand{\GL}{\mathrm{GL}}

\newcommand{\Gal}{\mathrm{Gal}}
\newcommand{\modulo}{\mathrm{mod}}
\newcommand{\Frob}{\mathrm{Frob}}

\newcommand{\Ker}{\mathrm{Ker}}

\newcommand{\cS}{\mathcal{S}}

\newcommand{\cO}{\mathcal{O}}

\newcommand{\ff}{\mathfrak{f}}

\newcommand{\fP}{\mathfrak{P}}

\newcommand{\Ind}{{\rm Ind}}
\newcommand{\PGL}{\mathrm{PGL}}

\newcommand{\CL}{\mathrm{CL}}

\newcommand{\pf}{{\bf Proof. }}
\newcommand{\qed}{\hspace* {.5cm} \hfill $\Box$}

\newcommand{\mat}[4]{
  \left( \begin{array}{cc} #1 & #2 \\ #3 & #4 \end{array} \right)}
\newcommand{\uomk}{\underline{\omega}^{\otimes k}}
\newcommand{\om}{{\underline{\omega}}}
\newcommand{\Om}{{\underline{\Omega}}}

\begin{document}

\title{Dihedral Galois representations\\
and Katz modular forms}
\author{Gabor Wiese\footnote{Mathematisch Instituut,
Universiteit Leiden,
Postbus 9512,
2300 RA Leiden,
The Netherlands, \newline
gabor@math.leidenuniv.nl,  {\tt http://www.math.leidenuniv.nl/$\sim$gabor/}}}
\maketitle

\begin{abstract}
We show that any two-dimensional odd dihedral representation~$\rho$ 
over a finite field of characteristic $p>0$
of the absolute Galois group of the rational numbers
can be obtained from a Katz modular form of
level~$N$, character~$\epsilon$ and weight~$k$, where $N$ is the
conductor, $\epsilon$ is the prime-to-$p$ part of the determinant and
$k$ is the so-called minimal weight of~$\rho$. In particular, $k=1$ if
and only if $\rho$ is unramified at~$p$.  Direct arguments are used in
the exceptional cases, where general results on weight and level
lowering are not available.

\end{abstract}

\section{Introduction}

In \cite{S} Serre conjectured that any odd irreducible continuous
Galois representation $\rho: G_\QQ \to \GL_2 (\Fpbar)$ for a prime $p$
comes from a modular form in characteristic~$p$ of a certain 
level $N_\rho$, weight $k_\rho \ge 2$ and character $\epsilon_\rho$.
Later Edixhoven discussed in
\cite{EW} a slightly modified definition of weight, the so-called 
{\em minimal weight}, denoted $k(\rho)$, by invoking Katz' theory of
modular forms.  In particular, one has that $k(\rho)=1$ 
if and only if $\rho$ is unramified at~$p$.

The present note contains a proof of this conjecture for
{\em dihedral representations}. 
We define those to be the continuous irreducible Galois representations 
that are induced from a character of the absolute Galois group of a 
quadratic number field. 
Let us mention that this is equivalent to imposing that the projective 
image is isomorphic to a dihedral group $D_n$ with $n \ge 3$.

\begin{thm}\label{dihthm}
Let $p$ be a prime and $\rho: G_\QQ \to \GL_2 (\Fpbar)$ an odd
dihedral representation.
As in \cite{S} define $N_\rho$ to be the conductor of~$\rho$ 
and $\epsilon_\rho$ to be the prime-to-$p$ part of $\det \circ \rho$
(considered as a character of $(\ZZ/(N_\rho p)\ZZ)^*$).
Define $k(\rho)$ as in~\cite{EW}. 

Then there exists a normalised Katz eigenform
$f \in \cS_{k(\rho)}(\Gamma_1(N_\rho), \epsilon_\rho, \Fpbar)_\Katz$, whose
associated Galois representation $\rho_f$
is isomorphic to~$\rho$.
\end{thm}

We will on the one hand show directly that $\rho$ comes from a 
Katz modular form of level~$N_\rho$, character~$\epsilon_\rho$ 
and minimal weight $k(\rho)=1$, if $\rho$ is unramified at~$p$.
If on the other hand $\rho$ is ramified at~$p$, we will
finish the proof by applying the fundamental work by Ribet, Edixhoven, Diamond, 
Buzzard and others on ``weight and level lowering''.

Let us recall that in weight at least~$2$ every Katz modular form 
on $\Gamma_1$ is classical, i.e.\ a reduction from a characteristic zero 
form of the same level and weight. 
Hence multiplying by the Hasse invariant, if necessary, it follows from 
theorem~\ref{dihthm} that every odd dihedral representation as above
also comes from a classical modular form of level~$N_\rho$ and
Serre's weight~$k_\rho$.
However, if one also wants the character to be~$\epsilon_\rho$, one
has to exclude in case $p=2$ that $\rho$ is induced 
from $\QQ(i)$ and in case $p=3$ that $\rho$ is induced from $\QQ(\sqrt{-3})$ 
(see \cite{B}, corollary~2.7, and \cite{D}, corollary~1.2).

Edixhoven's theorem on weight lowering (\cite{EW}, theorem~4.5) states
that modularity in level $N_\rho$ and the modified weight $k(\rho)$ 
follows from modularity in level $N_\rho$ and Serre's weight~$k_\rho$,
unless one is in a so-called {\em exceptional case}.
A representation $\rho: G_\QQ \to \GL_2(\Fpbar)$ is called {\em exceptional}
if the semi-simplification of its restriction to a decomposition group at~$p$ 
is the sum of two copies of an unramified character.
Because of work by Coleman and Voloch the only open case left is that
of characteristic~$2$ (see the introduction of~\cite{EW}).

Exceptionality at $2$ is a common phenomenon for mod~$2$ 
dihedral representations.
One can for example consider the quadratic field $K = \QQ(\sqrt{2089})$
and check that $2$ splits completely in the Hilbert class field of~$K$,
which has degree~$3$ over~$K$,
whence the associated dihedral Galois representation is exceptional.
Another type of example is provided by the dihedral representation
obtained from the Hilbert class field $H$ of $K = \QQ(\sqrt{229})$,
which also has degree~$3$.
The prime~$2$ stays inert in $\cO_K$, so $2 \cO_K$ splits completely
in~$H$, whence the restriction of the corresponding dihedral representation 
to a decomposition group at~$2$ is cyclic of order~$2$ and thus is 
also exceptional.

Let us point out that some of the weight one forms that we obtain
cannot be lifted to characteristic zero forms of weight one and the same 
level, so that the theory of modular forms by Katz becomes necessary. 
Namely, if $p=2$ and the dihedral representation in question has odd 
conductor~$N$ and is induced from a real quadratic field~$K$ of 
discriminant~$N$, 
whose fundamental units have norm~$-1$, then there does not exist an odd 
characteristic zero representation with conductor dividing~$N$
that reduces to~$\rho$. The representation 
coming from the quadratic field $\QQ(\sqrt{229})$ used above, can also here serve
as an example.

The fact that dihedral representations come from {\em some}
modular form is well-known (apparently already due to Hecke). So the subtle issue
is to adjust the level, character and weight.
It should be noted that Rohrlich and Tunnell solved many cases for $p=2$
with Serre's weight~$k_\rho$ by rather elementary means in~\cite{RT},
however, with the more restrictive definition of a dihedral representation
to be such that its image in $\GL_2(\Ftbar)$, and not in $\PGL_2(\Ftbar)$, 
is isomorphic to a dihedral group.

Let us also mention that it is possible to do computations of
weight one forms in positive characteristic on a computer (see~\cite{app})
and thus to collect evidence for Serre's conjecture in some cases.

This note is organised as follows.
The number theoretic ingredients on dihedral representations
are provided in section~2. 
In section~3 some results on oldforms, also in positive
characteristic, are collected.
Section~4 is devoted to the proof of theorem~\ref{dihthm}.
Finally, in section~5 we include a result on the irreducibility of certain
mod~$p$ representations.

\section{Dihedral representations}

We shall first recall some facts on Galois representations.
Let $\rho: G_\QQ \to \GL(V)$ be a continuous representation with
$V$ a $2$-dimensional vector space over an algebraically closed 
discrete field $k$.

Let $L$ be the number field such that $\Ker(\rho) = G_L$ (by the
notation $G_L$ we always mean the absolute Galois group of $L$).
Given a prime $\Lambda$ of $L$ dividing the rational prime $l$, we 
denote by $G_{\Lambda,i}$ the $i$-th ramification group in lower 
numbering of the local extension $L_\Lambda|\QQ_l$. Furthermore, one sets 
$$n_l(\rho) = \sum_{i\ge 0} \frac{\dim (V/ V^{G_{\Lambda,i}})}
   {(G_{\Lambda,0} : G_{\Lambda,i})}.$$
This number is an integer, which is independent of the choice of 
the prime $\Lambda$ above $l$.
With this one defines the {\em conductor} of $\rho$ to be
$\ff(\rho) = \prod_l l^{n_l(\rho)},$
where the product runs over all primes~$l$ different from the
characteristic of~$k$. If $k$ is the field of complex numbers, $\ff(\rho)$
coincides with the {\em Artin conductor}.

Let $\rho$ be a dihedral representation. Then $\rho$ is induced
from a character $\chi: G_K \to k^*$ for a quadratic number field~$K$
such that $\chi \neq \chi^\sigma$, with 
$\chi^\sigma (g) = \chi(\sigma^{-1}g\sigma)$ for all $g \in G_K$,
where $\sigma$ is a lift to $G_\QQ$ of the non-trivial element
of $G_{K|\QQ}$.
For a suitable choice of basis we then have the following explicit 
description of~$\rho$:
If an unramified prime~$l$ splits in~$K$ as $\Lambda \sigma(\Lambda)$,
then 
$ \rho(\Frob_l) = 
   \mat {\chi(\Frob_{\Lambda})} 0 0 {\chi^\sigma(\Frob_{\Lambda})}.$
Moreover, $\rho(\sigma)$~is represented by the matrix $\mat 0 1 {\chi(\sigma^2)} 0$.
As $\rho$ is continuous, its image is a finite group, say, of order~$m$.

\begin{lem}\label{dihlem}
Let $\rho: G_\QQ \to \GL_2(\Fpbar)$ be an odd dihedral representation
that is unramified at~$p$.
Define $K$, $\chi$, $\sigma$ and $m$ as above.
Let $N$ be the conductor of~$\rho$. 
Let $\zeta_m$ a primitive $m$-th root of unity 
and $\fP$ a prime of~$\QQ(\zeta_m)$ above~$p$.

Then one of the following two statements holds.
\begin{enumerate}[(a)]
\item There exists an odd dihedral representation
$\widehat{\rho}: G_\QQ \to \GL_2(\ZZ[\zeta_m])$, which has Artin conductor $N$
and reduces to $\rho$ modulo~$\fP$.
\item One has that $p=2$ and $K$ is real quadratic.
Moreover, there is an infinite set $S$ of primes such that for each $l \in S$
the trace of $\rho(\Frob_l)$ is zero, and there exists an odd dihedral representation
$\widehat{\rho}: G_\QQ \to \GL_2(\ZZ[\zeta_m])$,
which has Artin conductor $Nl$ and reduces to $\rho$ modulo~$\fP$.
\end{enumerate}
\end{lem}

\pf
Suppose that the quadratic field $K$ equals $\QQ(\sqrt{D})$ with
$D$ square-free.
The character $\chi: G_K \to k^*$ can be uniquely lifted to a character 
$\widetilde{\chi}: G_K \to \ZZ[\zeta_m]^*$ of the same order, which reduces
to $\chi$ modulo $\fP$. Denote by $\widetilde{\rho}$
the continuous representation $\Ind_{G_K}^{G_\QQ} \widetilde{\chi}$.
For the choice of basis discussed above the matrices representing $\rho$ can be 
lifted to matrices representing~$\widetilde{\rho}$, whose non-zero entries are 
in the $m$-th roots of unity.
Then for a subgroup $H$ of the image $\rho(G_\QQ)$, one has that
$(\Fpbar^2)^H$ is isomorphic to $(\ZZ[\zeta_m]^2)^H \otimes \Fpbar$.
Hence the conductor of~$\rho$ equals the Artin conductor 
of $\widetilde{\rho}$, as $\widetilde{\rho}$ is unramified at~$p$.
Alternatively, one can first remark that the conductor of $\chi$ equals
the conductor of $\widetilde{\chi}$ and then use the formulae
$\ff(\rho)=\Norm_{K|\QQ}(\ff(\chi)) D$ and 
$\ff(\widetilde{\rho})=\Norm_{K|\QQ}(\ff(\widetilde{\chi})) D$.

Thus condition (a) is satisfied if $\widetilde{\rho}$ is odd.
Let us now consider the case when $\widetilde{\rho}$ is even. This immediately
implies $p=2$ and that the quadratic field $K$ is real,
as is the number field $L$ whose absolute Galois group $G_L$ equals
the kernel of~$\rho$, and hence also the kernel of~$\widetilde{\chi}$.
We shall now adapt ``Serre's trick'' from \cite{RT}, p.~307, to
our situation.

Let $\ff$ be the conductor of~$\widetilde{\chi}$.
As $L$ is totally real, $\ff$ is a finite ideal of~$\cO_K$.
Via class field theory, $\widetilde{\chi}$ can be identified with 
a complex character of $\CL_K^\ff$, the ray class group modulo~$\ff$.
Let $\infty_1, \infty_2$ be the infinite places of $K$.
Consider the class
$$ c = [\{ (\lambda) \in \CL_K^{4D\ff\infty_1\infty_2} \; \mid \;
    \Norm(\lambda)<0, \lambda \equiv 1 \, \modulo \, 4D\ff\}]$$
in the ray class group of~$K$ modulo $4D\ff\infty_1\infty_2$.
By Cebotarev's density theorem the primes of $\cO_K$ are uniformly
distributed over the conjugacy classes of $\CL_K^{4D\ff\infty_1\infty_2}$.
Hence, there are infinitely many primes $\Lambda$ 
of degree~$1$ in the class~$c$.
Take $S$ to be the set of rational primes lying under them.
Let a prime $\Lambda$ from the class $c$ be given. 
It is principal, say $\Lambda = (\lambda)$, and coprime to $4D\ff$.
By construction we have 
$c^2 = [\Lambda^2] = 1$. As $\CL_K^\ff$ is a quotient of 
$\CL_K^{4D\ff\infty_1\infty_2}$, the class of $\Lambda$ in
$\CL_K^\ff$ has order $1$ or $2$. Since $p=2$, the character ${\chi}$
has odd order and we conclude that ${\chi}(\Lambda)=1$.

We have $\lambda \equiv 1 \, \modulo \, 4D\ff$ and $\Norm(\lambda) = -l$ 
for some odd prime~$l$.
Hence, the extension $K(\sqrt{\lambda})$ has signature $(2,1)$
and is unramified at~$2$ and at the primes dividing~$D\ff$.
We represent $K(\sqrt{\lambda})$ by the quadratic character
$\xi:G_K \to \{\pm 1\}$.
For the complex conjugation, the ``infinite Frobenius element'', $\Frob_{\infty_1}$,
we have that $\xi(\Frob_{\infty_1})\xi^\sigma(\Frob_{\infty_1}) =-1$.
We now consider the representation~$\widehat{\rho}$
obtained by induction from the character 
$\widehat{\chi}=\widetilde{\chi}\xi$.
Using the same basis as in the discussion at the beginning of this section,
an element $g$ of $G_K$ is represented by the matrix
$\mat {\widetilde{\chi}(g) \xi(g)} 0 0 {\widetilde{\chi}^\sigma(g) \xi^\sigma(g)}$.
In particular, we obtain that the determinant of 
$\Frob_\infty $ over $\QQ$ equals $-1$, whence $\widehat{\rho}$ is odd.
Moreover, as $l$ splits in~$K$, one has that $\rho(\Frob_l)$ is the identity
matrix, so that the trace of $\rho(\Frob_l)$ is zero.

The reduction of $\widehat{\rho}$ equals~$\rho$, as $\xi$ is trivial 
in characteristic~$2$.
Moreover, outside $\Lambda$ the conductor of $\widehat{\chi}$ equals
the conductor of~$\widetilde{\chi}$. At the prime~$\Lambda$ the
local conductor of~$\widehat{\chi}$ is~$\Lambda$, as the 
ramification is tame. Consequently, the Artin conductor of~$\widehat{\rho}$ 
equals~$Nl$.
\qed
\medskip

Also without the condition that it is unramified at~$p$, one can lift a 
dihedral representation to characteristic zero, however, losing
control of the Artin conductor.

\begin{lem}\label{dihlemeins}
Let $\rho: G_\QQ \to \GL_2(\Fpbar)$ be an odd dihedral representation.
Define $K$, $\chi$, $m$, $\zeta_m$ and $\fP$ as in the previous lemma.

There exists an odd dihedral representation
$\widehat{\rho}:G_\QQ \to \GL_2(\ZZ[\zeta_m])$, whose reduction modulo~$\fP$
is isomorphic to~$\rho$.
\end{lem}

\pf
We proceed as in the preceding lemma for the 
definitions of $\widetilde{\chi}$ and $\widetilde{\rho}$.
If $\widetilde{\rho}$ is even, then $p=2$ and $K$ is real.
In that case we choose some $\lambda \in \cO_K - \ZZ$,
which satisfies $\Norm(\lambda)<0$.
The field $K(\sqrt{\lambda})$ then has signature $(2,1)$ 
and gives a character $\xi: G_K \to \ZZ[\zeta_m]^*$.
As in the proof of the preceding lemma one obtains that the representation 
$\widehat{\rho} = \Ind_{G_K}^{G_\QQ} \widetilde{\chi}\xi$
is odd and reduces to~$\rho$ modulo~$\fP$.
\qed

\section{On oldforms}

In this section we collect some results on oldforms. We try to
stay as much as possible in the characteristic zero setting. However,
we also need a result on Katz modular forms.

\begin{prop}\label{oldzero}
Let $N,k,r$ be positive integers, $p$ a prime and $\epsilon$ 
a Dirichlet character of modulus~$N$.
The homomorphism
$$ \phi_{p^r}^N: \big(\cS_k(\Gamma_1(N),\epsilon,\CC)\big)^{r+1} \hookrightarrow
   \cS_k(\Gamma_1(Np^r),\epsilon,\CC), \;\;
   (f_0, f_1, \dots, f_r) \mapsto \sum_{i=0}^r f_i(q^{p^i})$$
is compatible with all Hecke operators $T_n$ with $(n,p)=1$.

Let $f \in \cS_k(\Gamma_1(N),\epsilon,\CC)$ be a normalised eigenform for 
all Hecke operators. Then the forms
$f(q), f(q^{p^2}), \dots, f(q^{p^r})$ in the image of $\phi_{p^r}^N$
are linearly independent, and on their span the action of the 
operator $T_p$ in level $Np^r$ is given by the matrix
$$ \left( \begin{array}{cccccc}  
a_p(f) & 1 & 0 & 0  & \dots   & 0 \\
-\delta p^{k-1}\epsilon(p) & 0 & 1 & 0 & \dots & 0 \\
0 & 0 & 0 & 1 & \dots & 0 \\
&&\vdots&&&\\
0 & \dots & 0 & 0 & 0 & 1 \\
0 & \dots & 0 & 0 & 0 & 0 \\
\end{array} \right),$$
where $\delta = 1$ if $p \nmid N$ and $\delta = 0$ otherwise.
\end{prop}

\pf
The embedding map and its compatibility with the Hecke action
away from $p$ is explained in \cite{DI}, section 6.1.
The linear independence can be checked on $q$-expansions.
Finally, the matrix can be elementarily computed.
\qed

\begin{cor}\label{cornull}
Let $p$ be a prime, $r \ge 0$ some integer and
$f \in \cS_k (\Gamma_1(Np^r), \epsilon, \CC)$ an eigenform
for all Hecke operators.
Then there exists an eigenform for all Hecke operators
$\tilde{f} \in \cS_k (\Gamma_1(Np^{r+2}), \epsilon, \CC)$,
which satisfies
$a_l(\tilde{f}) = a_l(f)$ for all primes $l \neq p$ and
$a_p(\tilde{f}) = 0$.
\end{cor}

\pf
One computes the characteristic polynomial of the operator
$T_p$ of proposition \ref{oldzero} and sees that it has $0$ as a root
if the dimension of the matrix is at least $3$. 
Hence one can choose the desired eigenform $\widetilde{f}$
in the image of~$\phi_{p^2}^{Np^r}$.
\qed
\medskip

As explained in the introduction, Katz' theory of modular forms 
ought to be used in the study of Serre's conjecture.
Following \cite{ES}, we briefly recall this concept, 
which was introduced by Katz in \cite{K}.
However, we shall use a ``non-compactified'' version.

Let $N\ge 1$ be an integer and $R$ a ring, in which $N$ is invertible.
One defines the category $[\Gamma_1(N)]_R$, whose objects are
pairs $(E/S/R,\alpha)$, where $S$ is an $R$-scheme, $E/S$ an elliptic
curve (i.e. a proper smooth morphism of $R$-schemes, whose geometric
fibres are connected smooth curves of genus one, together with
a section, the ``zero section'', $0: S \to E$) 
and $\alpha: (\ZZ/N\ZZ)_S \to E[N]$, the {\em level structure}, is an embedding of $S$-group
schemes. The morphisms in the category are cartesian diagrams
$$ \xymatrix@=.8cm{
E'  \ar@{->}[r] \ar@{}[rd]|{\Box}& E \\
S' \ar@{->}[r] \ar@{<-}[u]  & S, \ar@{<-}[u] } $$
which are compatible with the zero sections and the level structures.
For every such elliptic curve $E/S/R$ we let $\om_{E/S} = 0^*\Om_{E/S}$.
For every morphism $\pi: E'/S'/R \to E/S/R$
the induced map $\om_{E'/S'} \to \pi^* \om_{E/S}$ is an isomorphism.

A {\em Katz cusp form} $f \in \cS_k(\Gamma_1(N),R)_\Katz$ assigns to every
object $(E/S/R,\alpha)$ of $[\Gamma_1(N)]_R$ an element
$f(E/S/R,\alpha) \in \om_{E/S}^{\otimes k}(S)$, compatibly for the 
morphisms in the category, subject to the condition that all $q$-expansions
(which one obtains by adjoining all $N$-th roots of unity and
plugging in a suitable Tate curve) only have positive terms.

For the following definition let us remark that if $m \ge 1$ is coprime to $N$
and is invertible in $R$, then any morphism of group schemes of the form
$\phi_{Nm}: (\ZZ/Nm\ZZ)_S \to E[Nm]$
can be uniquely written as $\phi_N \times_S \phi_m$ with 
$\phi_N: (\ZZ/N\ZZ)_S \to E[N]$ and $\phi_m: (\ZZ/m\ZZ)_S \to E[m]$.

\begin{defi}
A Katz modular form $f \in \cS_k(\Gamma_1(Nm),R)_\Katz$ is
called {\em independent of $m$} if
for all elliptic curves $E/S/R$, all $\phi_N: (\ZZ/N)_S \hookrightarrow E[N]$
and all $\phi_m, \phi_m': (\ZZ/m)_S \hookrightarrow E[m]$ one has the
equality
$$ f(E/S/R,\phi_N \times_S \phi_m) = f(E/S/R,\phi_N \times_S \phi_m') 
  \in \uomk_{E/S}(S). $$
\end{defi}

\begin{prop}\label{propkatz}
Let $N$, $m$ be coprime positive integers and $R$ a ring, which contains
the $Nm$-th roots of unity and $\frac{1}{Nm}$.
A Katz modular form $f \in \cS_k(\Gamma_1(Nm),R)_\Katz$ is
independent of $m$ if and only if there exists a Katz modular form
$g \in \cS_k(\Gamma_1(N),R)_\Katz$ such that
$$ f(E/S/R, \phi_{Nm}) = g(E/S/R, \phi_{Nm} \circ \psi) $$
for all elliptic curves $E/S/R$ and all 
$\phi_{Nm}: (\ZZ/Nm\ZZ)_S \hookrightarrow E[Nm]$. Here $\psi$ denotes the canonical
embedding $(\ZZ/N\ZZ)_S \hookrightarrow (\ZZ/Nm\ZZ)_S$ of
$S$-group schemes.
In that case, $f$ and $g$ have the same $q$-expansion at~$\infty$.
\end{prop}

\pf
If $m = 1$, there is nothing to do.
If necessary replacing $m$ by~$m^2$, we can hence assume that
$m$ is greater equal~$3$.

Let us now consider the category $[\Gamma_1(N; m)]_R$, whose objects are triples
$(E/S/R,\phi_N,\psi_m)$, where $S$ is an $R$ scheme, $E/S$ an elliptic
curve, $\phi_N : (\ZZ/N\ZZ)_S \hookrightarrow E[N]$ an embedding of group
schemes and $\psi_m (\ZZ/m\ZZ)_S^2 \cong E[m]$ an isomorphism 
of group schemes. The morphisms are cartesian diagrams compatible with
the zero sections, the $\phi_N$ and the $\psi_m$ as before.

We can pull back the form $f \in \cS_k(\Gamma_1(Nm),R)_\Katz$ to a Katz
form $h$ on~$[\Gamma_1(N;m)]_R$ as follows. 
First let $\beta: (\ZZ/m\ZZ)_S \hookrightarrow (\ZZ/m\ZZ)_S^2$ be the embedding of
$S$-group schemes defined by mapping onto the first factor. 
Using this, $f$ gives rise to $h$ by setting
$$h((E/S/R,\phi_N,\psi_m)) = 
f((E/S/R,\phi_N,\psi_m \circ \beta)) \in \om_{E/S}^{\otimes k}(S).$$
As $f$ is independent of $m$, it is clear that $h$ is independent of $\psi_m$
and thus invariant under the natural $\GL_2(\ZZ/m\ZZ)$-action.

As $m \ge 3$, one knows that the category $[\Gamma_1(N; m)]_R$
has a final object $(E^\univ/Y_1(N;m)_R/R,\alpha^\univ)$.
In other words, $h$ is an $\GL_2(\ZZ/m\ZZ)$-invariant global section of
$\om_{E^\univ/Y_1(N;m)_R}^{\otimes k}$.
Since this $R$-module is equal to $\cS_k(\Gamma_1(N),R)_\Katz$
(see e.g.\ equation 1.2 of \cite{ES}, p.\ 210), we find some 
$g \in \cS_k(\Gamma_1(N),R)_\Katz$ such that
$ f(E/S/R, \phi_{Nm}) = g(E/S/R, \phi_{Nm} \circ \psi) $
for all $(E/S/R,\phi_{Nm})$.

Plugging in the Tate curve, one sees that the standard $q$-expansions of $f$
and $g$ coincide.
\qed

\begin{cor}\label{corkatz}
Let $N,m$ be coprime positive integers, $p$ a prime not dividing $Nm$ and
$\epsilon: (\ZZ/N\ZZ)^* \to \Fpbar$ a character.
Let $f \in \cS_k(\Gamma_1(Nm),\epsilon,\Fpbar)_\Katz$ be a Katz cuspidal
eigenform for all Hecke operators. 

If $f$ is independent of $m$, then there exists an eigenform for all Hecke operators
$g \in \cS_k(\Gamma_1(N),\epsilon,\Fpbar)_\Katz$ such that the associated
Galois representations $\rho_f$ and $\rho_g$ are isomorphic.
\end{cor}

\pf
From the preceding proposition we get a modular form 
$g \in \cS_k(\Gamma_1(N),\epsilon,\Fpbar)_\Katz$, noting that the character
is automatically good. Because of the
compatibility of the embedding map with the operators $T_l$ for primes
$l \nmid m$, we find that $g$ is an eigenform for these operators.
As the operators $T_l$ for primes $l \nmid m$ commute with the others,
we can choose a form of the desired type.
\qed

\section{Proof of the principal result}

We first cover the weight one case.

\begin{thm}\label{pthm}
Let $p$ be a prime and $\rho: G_\QQ \to \GL_2 (\Fpbar)$ an odd
dihedral representation of conductor~$N$, which is unramified at~$p$.
Let $\epsilon$ denote the character $\det \circ \rho$.

Then there exists a Katz eigenform $f$ in 
$\cS_1(\Gamma_1(N),\epsilon, \Fpbar)_\Katz$, 
whose associated Galois representation is isomorphic to~$\rho$.
\end{thm}

\pf
Assume first that part~(a) of lemma~\ref{dihlem} applies to~$\rho$, and let
$\widehat{\rho}$ be a lift provided by that lemma.
A theorem by Weil-Langlands (Theorem~1 of \cite{S1})
implies the existence of a newform $g$ in 
$\cS_1(\Gamma_1(N),\det \circ \widehat{\rho}, \CC)$,
whose associated Galois representation is isomorphic to~$\widehat{\rho}$.
Now reduction modulo a suitable prime above~$p$ yields the desired modular form.
In particular, one does not need Katz' theory in this case.

If part~(a) of lemma~\ref{dihlem} does not apply, then part~(b) does,
and we let $S$ be the infinite set of primes provided.
For each $l \in S$ the theorem of Weil-Langlands yields a newform
$f^{(l)}$ in $\cS_1 (\Gamma_1(Nl), \CC)$, whose associated Galois representation
reduces to $\rho$ modulo~$\fP$, where $\fP$ is the ideal from
the lemma.
Moreover, the congruence $a_q(f^{(l)}) \equiv 0 \, \modulo \, \fP$ holds for all primes
$q \in S$ different from~$l$.

From corollary \ref{cornull} we obtain Hecke eigenforms 
$\widetilde{f}^{(l)} \in \cS_1 (\Gamma_1(Nl^3), \CC)$
such that $a_l(\widetilde{f}^{(l)}) = 0$ and 
$a_q(\widetilde{f}^{(l)}) = a_q(f^{(l)}) \equiv 0 \, \modulo \, \fP$
for all primes $q \in S$, $q \neq l$.
Reducing modulo the prime ideal~$\fP$, we get eigenforms
$g^{(l)} \in \cS_1(\Gamma_1(Nl^3),\epsilon, \Fpbar)$,
whose associated Galois representations are isomorphic to~$\rho$.
One also has $a_q(g^{(l)})=0$ for all $q \in S$.

The coefficients $a_q(f^{(l)})$ for all primes $q \mid N$ appear in
the L-series of the complex representation $\rho_{f^{(l)}}$ associated
to~$f^{(l)}$. As the image of $\rho_{f^{(l)}}$
is isomorphic to a fixed finite group $G$, not depending on~$l$, there are only
finitely many possibilities for the value of $a_q(f^{(l)})$.
Hence the same holds for the $g^{(l)}$.
Consequently, there are two forms $g_1 = g^{(l_1)}$ and $g_2 = g^{(l_2)}$
for $l_1 \neq l_2$
that have the same coefficients at all primes $q \mid N$. 
For primes $q \nmid Nl_1l_2$ one has that the trace of
$\rho_{f^{(l_1)}}(\Frob_q)$ is congruent to the trace of 
$\rho_{f^{(l_2)}}(\Frob_q)$,
whence $a_q(g_1) = a_q(g_2)$. Let us point out that this includes
the case $q=p=2$, as the complex representation is unramified at~$p$.

In the next step we embed $g_1$ and $g_2$ into 
$\cS_1(\Gamma_1(Nl_1^3l_2^3),\epsilon, \Fpbar)_\Katz$
via the method in the statement of proposition \ref{propkatz}.
As the $q$-expansions coincide, $g_1$ and $g_2$ are mapped to the same form~$h$.
But as $h$ comes from $g_2$, it is independent of $l_1$
and analogously also of $l_2$.
Since $\rho_h = \rho$, theorem \ref{pthm} follows immediately from 
corollary \ref{corkatz}.
\qed
\medskip

We will deduce the cases of weight greater equal two from general results.
The current state of the art in ``level and weight lowering'' seems
to be the following theorem.

\begin{thm}\label{levellowering}[Ribet, Edixhoven, Diamond, Buzzard,$\dots$]
Let $p$ be a prime and $\rho: G_\QQ \to \GL_2(\Fpbar)$ a continuous
irreducible representation, which is assumed to come from some modular
form. Define $k_\rho$ and $N_\rho$ as in~\cite{S}.
If $p=2$, additionally assume either (i)~that the restriction of $\rho$
to a decomposition group at~$2$ is not contained within the
scalar matrices or (ii)~that $\rho$ is ramified at~$2$.

Then there exists a normalised eigenform 
$f \in \cS_{k_\rho}(\Gamma_1(N_\rho), \Fpbar)$ giving rise to $\rho$.
\end{thm}

\pf
The case $p \neq 2$ is theorem~1.1 of~\cite{D}, 
and the case $p=2$ with condition~(i)
follows from propositions 1.3 and~2.4 and theorem~3.2 of~\cite{B},
multiplying by the Hasse invariant if necessary.

We now show that if $p=2$ and $\rho$ restricted to a decomposition group $G_{\QQ_2}$
at~$2$ is contained within the scalar matrices, then $\rho$ is unramified at~$2$.
Let $\phi: G_\QQ \to \Ftbar^*$ be the character such that $\phi^2 = \det \circ \rho$.
As $\phi$ has odd order, it is unramified at~$2$ because of the Kronecker-Weber
theorem. If $\rho$ restricted to $G_{\QQ_2}$
is contained within the scalar matrices, then we have that $\rho|_{G_{\QQ_2}}$
is $\mat {\phi|_{G_{\QQ_2}}} 0 0 {\phi|_{G_{\QQ_2}}}$, whence $\rho$ is
unramified at~$2$.
\qed
\medskip

{\bf Proof of theorem \ref{dihthm}.}
Let $\rho$ be the dihedral representation from the assertion.
If $\rho$ is unramified at~$p$, one has $k(\rho)=1$, and theorem~\ref{dihthm}
follows from theorem~\ref{pthm}.

If $\rho$ is ramified at~$p$, then let $\widehat{\rho}$ be a characteristic zero 
representation lifting~$\rho$,
as provided by lemma~\ref{dihlemeins}. The theorem by Weil-Langlands already used above
(Theorem 1 of \cite{S1})
implies the existence of a newform in weight one and characteristic 
zero giving rise to $\widehat{\rho}$.
So from theorem~\ref{levellowering} we obtain that $\rho$ comes from a modular
form of Serre's weight~$k_\rho$ and level~$N_\rho$.
Let us note that using Katz modular forms the character is automatically 
the conjectured one~$\epsilon_\rho$.

The weights $k_\rho$ and $k(\rho)$ only differ in two cases
(see \cite{EW}, remark 4.4). The first case is when $k(\rho)=1$.
The other case is when $p=2$ and $\rho$ is not finite at~$2$.
Then one has $k(\rho) = 3$ and $k_\rho = 4$. In that case one
applies theorem~3.4 of~\cite{EW} to obtain an eigenform
of the same level and character in weight~$3$, or one applies
theorem~3.2 of~\cite{B} directly.
\qed

\section{An irreducibility result}

We first study the relation between the level of an eigenform in 
characteristic~$p$ and the conductor of the associated Galois 
representation.

\begin{lem}\label{lembound}
Let $\rho: G_\QQ \to \GL_2(\Fpbar)$ be a continuous representation of
conductor~$N$, and let $k$~be a positive integer.
If $f \in \cS_k(\Gamma_1(M),\epsilon,\Fpbar)_\Katz$ is a Hecke eigenform
giving rise to $\rho$, then $N$ divides $M$.
\end{lem}

\pf
By multiplying with the Hasse invariant, if necessary, we can
assume that the weight is greater equal~$2$. Hence the form $f$
can be lifted to characteristic zero (see e.g.\ \cite{DI}, Theorem~12.3.2)
in the same level.
Thus there exists a newform $g$, say of level $L$, 
whose Galois representation $\rho_g$ reduces to~$\rho$.
Now Proposition~0.1 of \cite{L} yields that $N$ divides~$L$.
As $L$ divides $M$, the lemma follows.
\qed
\medskip

We can derive the following proposition, which is of independent interest.

\begin{prop}\label{propirr}
Let $f \in \cS_k(\Gamma_0(N),\Fpbar)_\Katz$ be a normalised Hecke
eigenform for a square-free level $N$ with $p \nmid N$ in some 
weight~$k \ge 1$.
\begin{enumerate}[(a)]
\item If $p=2$, the associated Galois representation is either 
irreducible or trivial.
\item For any prime $p$ the associated Galois representation is either 
irreducible or corresponds to the sum of a character of $\QQ(\zeta_p)$
and its inverse, where $\zeta_p$ is a primitive $p$-th root of unity.
\end{enumerate}
\end{prop}

\pf
Let us assume that the representation $\rho$ associated to $f$ 
is reducible. Since $\rho$ is semi-simple, 
it is isomorphic to the direct sum of two characters $\chi_1 \oplus \chi_2$.
As the determinant is $1$, we have that $\chi_2^{-1} = \chi_1$.
Consequently, the conductor of $\rho$ is the square of the conductor
of~$\chi_1$. Lemma \ref{lembound} implies that the 
conductor of~$\rho$ divides~$N$. As we have assumed this
number to be square-free, we have that $\rho$ can only ramify in~$p$.

The number field $L$ with $G_L = \Ker(\rho) = \Ker(\chi_1)$ is cyclic.
As only $p$ can be ramified, it follows that $L$ is contained in 
$\QQ(\zeta_{p^n})$ for some $p^n$-th root of unity.
Since the order of $\chi_1$ is prime to~$p$, we conclude that
$L$ is contained in $\QQ(\zeta_p)$.

So the trace of Frobenius is the sum of a $(p-1)$-th root of unity
and its inverse. In characteristic~$2$ this implies that all traces
are zero, whence $\rho$ is the trivial representation.
\qed

\end{document}